\documentclass[a4paper,12pt]{amsart}

\usepackage[utf8]{inputenc}
\usepackage{graphicx}

\usepackage{multirow}
\usepackage{amssymb}
\usepackage{amsmath}
\usepackage{amsthm}

\usepackage[colorlinks]{hyperref}
\usepackage{hyperref}
\renewcommand\eqref[1]{(\ref{#1})} 

\allowdisplaybreaks \numberwithin{equation}{section}

\theoremstyle{plain}
\newtheorem{theorem}{Theorem}[section]
\newtheorem{prop}[theorem]{Proposition}
\newtheorem{problem}[theorem]{Problem}

\theoremstyle{definition}
\newtheorem{defn}[theorem]{Definition}
\newtheorem{rem}[theorem]{Remark}

\begin{document}

\title[Inverse source problem generated by the Dunkl operator]{One inverse source problem generated by the Dunkl operator}

\author[B. Bekbolat ]{Bayan Bekbolat}
\address{
  Bayan Bekbolat :
  \endgraf
  Al-Farabi Kazakh National University
  \endgraf
  Almaty, Kazakhstan
  \endgraf
  and
  \endgraf 
  Department of Mathematics: Analysis, Logic and Discrete Mathematics
  \endgraf
  Ghent University, Belgium
  \endgraf
  and
  \endgraf   
  Institute of Mathematics and Mathematical Modeling
  \endgraf
  Almaty, Kazakhstan
  \endgraf  
  and
  \endgraf 
  Suleyman Demirel University
  \endgraf
  Kaskelen, Kazakhstan
  \endgraf
  {\it E-mail address:} {\rm bekbolat@math.kz}
  }

\author[N. Tokmagambetov ]{Niyaz Tokmagambetov }
\address{
  Niyaz Tokmagambetov:
  \endgraf 
  Centre de Recerca Matem\'atica
  \endgraf
  Edifici C, Campus Bellaterra, 08193 Bellaterra (Barcelona), Spain
  \endgraf
  and
  \endgraf   
  Institute of Mathematics and Mathematical Modeling
  \endgraf
  125 Pushkin str., 050010 Almaty, Kazakhstan
  \endgraf  
  {\it E-mail address:} {\rm tokmagambetov@crm.cat; tokmagambetov@math.kz}
  }

\date{\today}

\thanks{This research was funded by the Science Committee of the Ministry of Science and Higher Education of the 
Republic of Kazakhstan (Grant No. AP14972634).}


\maketitle

\begin{abstract}
The aim of this paper is to study time-fractional pseudo-parabolic type equations generated by the Dunkl operator.
The forward problem is considered and its well-posedness is established. In particular, a prior estimates are obtained in the Sobolev type spaces and, explicit formulas for solutions of the problems are derived. Here we also deal with the left-sided Caputo fractional time derivative. 

As an application, we investigate an inverse source problem. Existence and uniqueness of the solution is proved. Moreover, we show that a solution pair is continuously depending on the initial and additional data, finalizing with a numerical test.
\end{abstract}

\section{Introduction}

In this paper we are interested in studying direct and inverse source problems (ISPs) for the time-fractional pseudo-parabolic equation of the form 
\begin{equation*}
\mathbb{D}_{0^+,t}^\gamma\left( u(t,x)-a\Lambda_{\alpha,x}^2 u(t,x) \right) -\Lambda_{\alpha,x}^2 u(t,x) + mu(t,x)  = f(x), \quad (t,x)\in Q_T,
\end{equation*}
generated by the Dunkl operator $\Lambda_{\alpha}$, where $Q_T=\{(t,x) : 0<t<T<+\infty, x\in\mathbb{R}\}$, $a,m>0$, and
\begin{equation*}
\mathbb{D}_{0^+,t}^\gamma =
\begin{cases}
     \mathcal{D}_{0^+,t}^\gamma, \quad 0<\gamma<1,\\
     \frac{d}{dt}, \quad \gamma=1,
\end{cases}
\end{equation*}
$\mathcal{D}_{0^+,t}^\gamma, 0<\gamma<1,$ is the left-sided Caputo fractional derivative.

The main scope of this paper is to investigate ISPs. In 1980, Rundell and Colton firstly studied ISP in \cite{Run80}. They considered the  evolution type equation
\begin{equation} \label{DiffEqu}
    \frac{du}{dt}+Au=f
\end{equation}
in a Banach space $X$, where $A$ is a linear operator in $X$ and $f$ is a constant vector in $X$, with conditions
\begin{equation*}
    u(0)=u_0, \quad \text{and} \quad u(T)=u_1.
\end{equation*}
Rundell and Colton proved a general theorem about the existence of a unique solution pair $(u(t),f)$ of the problem, using semigroups of operators. Then the results was applied to equations of parabolic and pseudo-parabolic types. In 1995 Bushuyev in \cite{B95}  considered the problem with an unknown source that can depend on time. However, in general case the problem was not realistic, so he considered the problem under a sufficient condition. The problem statement is 
\begin{equation*}
    \frac{\partial u}{\partial t} + Au = F(x,t) \quad \text{on} \quad \Omega\times(0,T)
\end{equation*}
\begin{equation*}
    D_v^j u = 0, \quad j\leq m-1 \quad \text{on} \quad \partial\Omega\times(0,T)
\end{equation*}
\begin{equation*}
    u = 0, \quad \text{on} \quad \Omega\times\{0\}, \Omega\times\{T\},
\end{equation*}
where $A$ is the linear elliptic partial differential operator of order $2m$ with the bounded measurable coefficients such that
\begin{equation*}
    (A\varphi,\varphi)\geq \|\varphi\|^2 
\end{equation*}
for all $\varphi\in H^{2m}(\Omega)\cap H_0^{m}(\Omega)$, $\mu=constant>0$. In 2004 in \cite{YG03}, Yaman and G\"{o}z\"{u}kızıl studied asymptotic behaviour of the solution of the inverse source problem for the pseudo-parabolic equation
\begin{equation*}
    (u(x,t) - \Delta u(x,t))_t - \Delta u(x,t) + \alpha u(x,t) = f(t)g(x,t), \quad Q_\infty = \Omega\times(0,\infty)
\end{equation*}
\begin{equation*}
    u(x,0) = u_0(x), x\in\Omega,
\end{equation*} 
\begin{equation*}
    u(x,t) = 0, \quad (x,t)\in\partial\Omega\times[0,\infty)
\end{equation*}
\begin{equation*}
    \int_\Omega u(x,t)(w-\Delta w)(x)dx=\varphi(t), \quad t\in[0,\infty),
\end{equation*}
where $g(x,t), w(x), \varphi(x), u_0(x)$, and the constant $\alpha$ are given, while $u(x,t)$ and $f(t)$ are unknown functions.

In 2011, K. Sakamoto and M. Yamamoto in \cite{SY11} considered an initial-boundary problem for a time fractional parabolic equation:
\begin{equation*}
    {}^c D_t^\alpha u(x,t) = r^\alpha (Lu)(x,t) + f(x)h(x,t), \quad x\in\Omega,\,\ t\in(0,T),\,\ 0<\alpha<1,
\end{equation*}
\begin{equation*}
    u(x,0) = 0, \quad x\in\Omega,
\end{equation*}
\begin{equation*}
    u(x,t) = 0, \quad x\in\partial\Omega, \,\ t\in(0,T),
\end{equation*}
where ${}^c D_t^\alpha, 0<\alpha<1$, is the Caputo derivative defined by
\begin{equation*}
    {}^c D_t^\alpha g(t) = \frac{1}{\Gamma(1-\alpha)} \int_0^t (t-\tau)^{-\alpha} \frac{d}{d\tau}g(\tau)d\tau
\end{equation*}
and $L$ is a symmetric uniformly elliptic operator
\begin{equation*}
    Lu = \sum_{i,j=1}^d \frac{\partial}{\partial x_i}\left(a_{ij}(x)\frac{\partial u}{\partial x_j}(x)\right) + p(x)u, \quad x\in\Omega.
\end{equation*}
The equation is commonly referred to as a fractional diffusion equation and serves as a model equation for describing diffusion within a particular class of porous media. The authors proved that the inverse problem, final overdetermining data, is well-posed in the Hadamard sense except for a discrete set of values of diffusion constants. In \cite{S13}, Slodi\u{c}ka considered ISP for the equation \eqref{DiffEqu} defined by
\begin{equation*}
    \frac{\partial u}{\partial t} + Au = f, \quad \text{in} \quad \Omega\times(0,T),
\end{equation*}
\begin{equation*}
    u = 0, \quad \text{on} \quad \Gamma\times(0,T),
\end{equation*}
\begin{equation*}
    u(x,0) = u_0(x), \quad \text{for} \quad x\in\Omega,
\end{equation*}
where $\Omega$ is a non-homogeneous and non-isotropic body, $\Gamma$ is a sufficiently smooth boundary, $T>0$ is the final time, and $A$ is a linear differential operator of second-order, strongly elliptic, and the right-hand side $f$ is assumed to be separable in both variables $x$ and $t$, i.e. $f(x,t)=g(x)h(t)$ (in this case $h(t)$ is unknown). The primary aim of the research was to identification the unknown function $h(t)$ based on additional data, under the assumption that $g(x)$ is already known.

We address the reader to a few articles (but not all) on the solvability of the inverse problems for the diffusion and sub-diffusion equations (\cite{CNYY09, JR15, KS10, KST17, OS12a, OS12b,  RTT19}) and fractional diffusion equations (\cite{SSB19, TT17, WYH13}).

Solvability of an ISP for pseudo-parabolic equation
\begin{equation*}
    \mathcal{D}_t^\alpha(u(t)+\mathcal{L}u(t))+\mathcal{M}u(t)=f(t) \quad \text{in} \quad \mathcal{H},
\end{equation*}
with Caputo fractional derivative $\mathcal{D}_t^\alpha$, of order $0<\alpha\leq1$, is considered in \cite{RSTT22} in 2022, where operators $\mathcal{L}$ and $\mathcal{M}$ are operators with the corresponding discrete spectra $\{\lambda_\xi\}_{\xi\in\mathcal{I}}$ and $\{\mu_\xi\}_{\xi\in\mathcal{I}}$ on $\mathcal{H}$. In \cite{RSTT22} well-posedness of direct and inverse
problems are obtained. The difference between our work and this particular work is that we consider an ISP for the pseudo-parabolic equation associated with the Dunkl operator, which is an operator with continuous spectrum. Therefore, our analysis is based on results of the Dunkl analysis. One dimensional Dunkl operator is the differential-difference operator given by the expression
\begin{equation*}
\Lambda_\alpha f(x)=\frac{d}{d x}f(x) + \left(\alpha+\frac{1}{2}\right)\frac{f(x)-f(-x)}{x},
\end{equation*}
where $x\in\mathbb{R}$ and $\alpha\geq-1/2$. This operator was introduced in 1989 by Dunkl \cite{D89}. It is associated with the reflection group $\mathbb{Z}_2$ on $\mathbb{R}$. These operators are very important in pure mathematics and provide a useful tool in the study of special functions associated with root systems. General references are \cite{D88, D89, D90, D91, D92, Jeu93, O93, R03} and \cite{A17}.


Let us give a brief description of our work. We are interested in studying the ISP for the time-fractional pseudo-parabolic equation   
\begin{equation}\label{EQ:Heat Dunkl}
\mathbb{D}_{0^+,t}^\gamma\left( u(t,x)-a\Lambda_{\alpha,x}^2 u(t,x) \right) -\Lambda_{\alpha,x}^2 u(t,x) + mu(t,x)  = f(x), \quad (t,x)\in Q_T,
\end{equation}
associated with the Dunkl operator $\Lambda_{\alpha}$, where $Q_T=\{(t,x) : 0<t<T<+\infty, x\in\mathbb{R}\}$, $a,m>0$, and
\begin{equation*}
\mathbb{D}_{0^+,t}^\gamma =
\begin{cases}
     \mathcal{D}_{0^+,t}^\gamma, \quad 0<\gamma<1,\\
     \frac{d}{dt}, \quad \gamma=1,
\end{cases}
\end{equation*}
$\mathcal{D}_{0^+,t}^\gamma, 0<\gamma<1,$ is the left-sided Caputo fractional derivative. We prove that the ISP for \eqref{EQ:Heat Dunkl} has a unique solution, pair of functions $(u,f)$, and the solution depends continuously on data and parameters, well-posedness in the sense of Hadamard. It is also shown that the direct problem to \eqref{EQ:Heat Dunkl} has a unique solution, where right hand side of the equation depends also from $t$ variable. In the case $a=0$, the time-fractional pseudo-parabolic equation \eqref{EQ:Heat Dunkl} reduces to the time-fractional heat equation
\begin{equation} \label{Rheat equation}
    \mathbb{D}_{0^+,t}^\gamma u(t,x) -\Lambda_{\alpha,x}^2 u(t,x) + mu(t,x)  = f(x), \quad (t,x)\in Q_T.
\end{equation}
In \cite{BST22}, B. Bekbolat, D. Serikbaev and N. Tokmagambetov considered the ISP for the equation \eqref{Rheat equation}. Authors in \cite{BST22} obtained the well-posedness results and a prior estimates in the Sobolev type spaces for the direct and inverse source problem for the equation \eqref{Rheat equation}. So, in this paper we are interested considering only case when $a>0$.

This paper is organized as follows. In Sect. \ref{Sect.Pre}, we give a brief exposition of analysis associated with Dunkl operator on $\mathbb{R}$. In Sect. \ref{Sect.Dir}, it is shown that the direct problem for the time-fractional pseudo-parabolic equation associated with Dunkl operator on $\mathbb{R}$ has a unique solution. In Sect. \ref{Sect.In}, our main results are stated and proved: the inverse source problem for the time-fractional pseudo-parabolic equation associated with Dunkl operator on $\mathbb{R}$ has a unique solution, stability analysis for inverse source problem and example for inverse source problem. 

Throughout this paper we use standard notations. The spaces $C^\infty(\mathbb{R})$, $\mathcal{S}(\mathbb{R})$, and $\mathcal{S}'(\mathbb{R})$ are the space of the infinitely differentiable functions on $\mathbb{R}$, the space of the Schwartz functions on $\mathbb{R}$, the space of the tempered distributions on $\mathbb{R}$ as the topological dual of $\mathcal{S}(\mathbb{R})$, respectively.

\section{Preliminaries}
\label{Sect.Pre}

\subsection{Dunkl operator}
(\cite[p.19]{A17}) For a real parameter $\alpha \geq -1/2$, we consider the differential-difference operator associated with the reflexion group $\mathbb{Z}_2$ on $\mathbb{R}$ defined by
\begin{equation}\label{EQ: Dunkl}
\Lambda_\alpha f(x)=\frac{d}{d x}f(x) + \left(\alpha+\frac{1}{2}\right)\frac{f(x)-f(-x)}{x}, \quad f\in C^1(\mathbb{R}),
\end{equation}
called Dunkl operator. If $\alpha=-1/2$, it turns into the ordinary differential operator $\Lambda_{-1/2} = d/dx$.

\begin{rem} (\cite[p.19]{A17}) 
From \eqref{EQ: Dunkl}, we can calculate
\begin{equation*}
    \Lambda_\alpha^2 f(x) = \frac{d^2}{d x^2}f(x) + \frac{2\alpha+1}{x}\cdot\frac{d}{d x}f(x) - \left(\alpha+\frac{1}{2}\right)\frac{f(x)-f(-x)}{x^2},
\end{equation*}
for all $f\in C^2(\mathbb{R})$. So, if $\alpha=-1/2$, we obtain a second order differential operator $\Lambda_{-1/2}^2 = d^2/dx^2$.
\end{rem}

\begin{prop} \cite[Proposition 3.5, p.19]{A17}
The Dunkl operators map the function spaces $C^\infty(\mathbb{R})$, $C_c^\infty(\mathbb{R})$, and $\mathcal{S}(\mathbb{R})$ into themselves. 
\end{prop}

For $\lambda\in\mathbb{R}$, the spectral problem
\begin{equation*}
    \Lambda_\alpha f(x)=(i\lambda)f(x), \quad f(0)=1, \quad x\in\mathbb{R},
\end{equation*}
has a unique solution $x\mapsto D_{\alpha}(ix\lambda)$ called Dunkl kernel given by
\begin{equation}\label{EQ: self lambda solution}
D_{\alpha}(i x\lambda)=j_\alpha(ix\lambda)+\frac{ix\lambda}{2(\alpha+1)}j_{\alpha+1}(ix\lambda), \quad x\in\mathbb{R},
\end{equation}
where
\begin{equation*}
j_\alpha(ix\lambda)=\sum_{k=0}^{\infty}\frac{\Gamma(\alpha+1)}{k!\Gamma(k+\alpha+1)}\left(\frac{ix\lambda}{2}\right)^{2k}
\end{equation*}
is the normalized Bessel function of order $\alpha$ \cite[Example 2.1, p.107]{R03}.
\begin{rem}
In the case $\alpha=-\frac{1}{2}$, we have the spectral problem $\frac{d}{dx}u(x)=(i\lambda)u(x)$, $u(0) = 1$, which has the (unique) solution $e^{ix\lambda}=D_{-1/2}(ix\lambda)$.
\end{rem}
\begin{defn}
We denote by $ L^p(\mathbb{R}, \mu_\alpha), 1\leq p\leq +\infty $, the space of measurable functions $f$ on $\mathbb{R}$ such that
\begin{equation*}
\|f\|_{p,\alpha}=\left( \int_\mathbb{R}|f(x)|^p d{\mu_\alpha(x)} \right)^{\frac{1}{p}}<+\infty, \quad 1\leq p<+\infty,
\end{equation*}
and
\begin{equation*}
\|f\|_{\infty}=\text{ess} \sup_{x\in\mathbb{R}}|f(x)|<+\infty,
\end{equation*}
where $\mu_\alpha$ is a measure defined on $\mathbb{R}$ by
\begin{equation*}
d\mu_\alpha(x)=\frac{|x|^{2\alpha+1}}{2^{\alpha+1}\Gamma(\alpha+1)}d x,\quad \alpha\geq-\frac{1}{2}.    
\end{equation*}
\end{defn}
For  $f\in L^1(\mathbb{R}, \mu_\alpha)$ the Dunkl transform is defined by
\begin{equation}
\label{Dunkl transform}
\mathcal{F}_{\alpha}(f)(\lambda) = \widehat{f}(\lambda) :=\int_\mathbb{R} f(x)D_{\alpha}(-i x\lambda) d\mu_\alpha(x), \quad \lambda\in\mathbb{R}.
\end{equation}
This transform has the following properties (\cite{Jeu93}):\\
i) For all $f\in L^1(\mathbb{R}, \mu_\alpha)$, Dunkl transform $\mathcal{F}_{\alpha}$ is a continuous function on $\mathbb{R}$ satisfying
\begin{equation*}
\|\mathcal{F}_{\alpha}(f)\|_{\infty}\leq\|f\|_{1,\alpha}.
\end{equation*}
ii) For all $f\in\mathcal{S}(\mathbb{R})$, we have
\begin{equation}\label{Commutation formula}
\mathcal{F}_{\alpha}(\Lambda_\alpha f)(\lambda) = i\lambda \mathcal{F}_{\alpha}(f)(\lambda), \quad \lambda\in\mathbb{R}.
\end{equation}
iii) ($L^1$-inversion) Let $f$ be a function in $L^1(\mathbb{R}, \mu_\alpha)$ such that $\mathcal{F}_{\alpha}(f)\in L^1(\mathbb{R}, \mu_\alpha)$. Then we have
\begin{equation}
\label{Invers Dunkl Transform}
    \mathcal{F}_{\alpha}^{-1}(f)(x) = \mathcal{F}_{\alpha}(f)(-x) \quad \text{a.e.} \quad  x\in\mathbb{R}.
\end{equation}
Moreover  
\begin{equation*}
    f(x) = \int_\mathbb{R} \mathcal{F}_{\alpha}(f)(\lambda)D_{\alpha}(ix\lambda) d\mu_{\alpha}(\lambda).
\end{equation*}
iv) (Plancherel theorem) Dunkl transform $\mathcal{F}_{\alpha}$ is an isometric isomorphism of $L^2(\mathbb{R}, \mu_{\alpha})$. In particular,
\begin{equation} \label{PTh}
    \|\mathcal{F}_{\alpha}(f)\|_{2,\alpha} = \|f\|_{2,\alpha}.
\end{equation}
\begin{defn}(\cite{MT04,S04}) For $s\in\mathbb{R}$, we define the space $\mathcal{H}_\alpha^2(\mathbb{R},\mu_\alpha)$ as the set of distributions $f\in\mathcal{S}'(\mathbb{R})$ satisfying
\begin{equation*}
    \int_\mathbb{R}|(1+\lambda^2)\mathcal{F}_{\alpha}(f)(\lambda)|^2d\mu_{\alpha}(\lambda)<\infty,
\end{equation*}
and norm of this space is given by the formula
\begin{equation*}
    \|f\|_{\mathcal{H}_\alpha^2}^2 := \int_\mathbb{R}|(1+\lambda^2)\mathcal{F}_{\alpha}(f)(\lambda)|^2d\mu_{\alpha}(\lambda).
\end{equation*}
\end{defn}
\begin{defn}
Let $u(t,\cdot)\in \mathcal{H}_\alpha^2(\mathbb{R},\mu_{\alpha})$. We denote by $C([0,T], \mathcal{H}_\alpha^2(\mathbb{R},\mu_{\alpha}))$ the space of continuous functions $u(\cdot, x)$ on [0,T] with the norm
\begin{equation*}
    \|u\|_{C([0,T], \mathcal{H}_\alpha^2(\mathbb{R},\mu_{\alpha}))} := \max_{0<t<T}\|u(t,\cdot)\|_{\mathcal{H}_\alpha^2(\mathbb{R},\mu_{\alpha})} < +\infty.
\end{equation*}
\end{defn}

\begin{defn} \cite[p. 18, Definition 3]{CF18}
\label{Defn}
Let $X$ be a Banach space. We say that $u\in C^\gamma([0,T],X)$ if $u\in C([0,T],X)$ and $\mathcal{D}_t^\gamma u\in C([0,T],X)$.
\end{defn}

\section{Direct problem}
\label{Sect.Dir}

This section deals with the Cauchy problem for the equation \eqref{EQ:Heat Dunkl}. The
Cauchy problem for the heat equation associated with the Dunkl operator 
\begin{equation*}
    \begin{cases}
    \Lambda_{\alpha,x}^2 u(t,x) - u_t(t,x) = 0\\
    u(0, t) = g(x)
    \end{cases}
\end{equation*}
was considered by Margit R\"{o}sler \cite{R98} on a domain $(0,\infty)\times \mathbb{R}$ with initial data $g\in C_b(\mathbb{R})$. Then nonhomogeneous problem
\begin{equation*}
    \begin{cases}
    u_t(t,x) - \Lambda_{\alpha,x}^2 u(t,x) = f(t,x)\\
    u(0, t) = g(x)
    \end{cases}
\end{equation*}
was considered by Hatem Mejjaoli \cite{M12,M13} on a domain $(0,\infty)\times \mathbb{R}$ when g belongs to homogeneous and nonhomogeneous Dunkl–Besov spaces.

\begin{problem} \label{Problem 1}
Let $0<\gamma\leq1$. Our aim is to find the function $u$ satisfying the equation 
\begin{equation}
\label{FracEq1}
    \mathbb{D}_{0^+,t}^\gamma \left(u(t,x) - a\Lambda_{\alpha,x}^2 u(t,x)\right) - \Lambda_{\alpha,x}^2u(t,x) + mu(t,x)  = f(t,x), \quad (t,x)\in Q_T,
\end{equation}
under the initial condition
\begin{equation}
\label{FracCon1}
    u(0,x)=g(x), \quad x\in\mathbb{R},
\end{equation}
where $f$ and $g$ are sufficiently smooth functions.
\end{problem}

The following theorem shows that Problem \ref{Problem 1} has a unique solution in the space $ C^\gamma([0,T],L^2(\mathbb{R},\mu_\alpha))\cap C([0,T],\mathcal{H}_\alpha^2(\mathbb{R},\mu_\alpha))$, where $0<\gamma\leq1$. 

\begin{theorem} \label{FracDirectTheorem}
a) Let $0<\gamma<1$. Assume that $f\in C^1([0,T],L^2(\mathbb{R},\mu_\alpha))$ and $g\in \mathcal{H}_\alpha^2(\mathbb{R},\mu_\alpha)$. Then the solution of the Problem \ref{Problem 1} exists, is unique, and given by the expression
\begin{equation} \label{expression}
\begin{split}
    u(t,x) &= \int_\mathbb{R}\int_\mathbb{R} \mathbb{E}_{\gamma,1}\left(-\frac{m+\lambda^2}{1+a\lambda^2}t^\gamma\right)g(y)D_\alpha(ix\lambda)D_\alpha(-iy\lambda)d\mu_\alpha(y)d\mu_\alpha(\lambda)\\
    &+ \int_\mathbb{R}\int_\mathbb{R}\int_0^t  (t-\tau)^{\gamma-1}\mathbb{E}_{\gamma,\gamma}\left(-\frac{m+\lambda^2}{1+a\lambda^2}(t-\tau)^\gamma\right)\frac{f(\tau,y)}{1+a\lambda^2}\\
    &\times D_\alpha(ix\lambda)D_\alpha(-iy\lambda)d\tau d\mu_\alpha(y)d\mu_\alpha(\lambda)
\end{split}
\end{equation}
where $\mathbb{E}_{\gamma,1}(t)$ is the classical Mittag-Leffler function and $\mathbb{E}_{\gamma,\gamma}(t)$ is the Mittag-Leffler type function:
\begin{equation*}
    \mathbb{E}_{\gamma,1}(t) := \sum_{k=0}^\infty \frac{t^k}{\Gamma(\gamma k+1)} \quad \mathbb{E}_{\gamma,\gamma}(t) := \sum_{k=0}^\infty \frac{t^k}{\Gamma(\gamma k+\gamma)}.
\end{equation*}
If $\gamma=1$, then $\mathbb{E}_{1,1}(t)=e^t$.\\
b) Let $\gamma=1$. Assume that $f\in C([0,T],L^2(\mathbb{R},\mu_{\alpha}))$ and $g\in \mathcal{H}_\alpha^2(\mathbb{R},\mu_{\alpha})$. Then the Problem \ref{Problem 1} has a unique solution, which is given by the expression \eqref{expression}.
\end{theorem}

For more information about the classical Mittag-Leffler function $\mathbb{E}_{\gamma,1}(t)$ and the Mittag-Leffler type function $\mathbb{E}_{\gamma,\gamma}(t)$ see e.g. \cite[p. 40 and p. 42]{KST06}.

In \cite{Sim14} the following estimate for the Mittag-Leffler function is proved, when $0<\gamma<1$ (not true for $\gamma\geq1$)
\begin{equation*}
    \frac{1}{1+\Gamma(1-\gamma)t}\leq\mathbb{E}_{\gamma,1}(-t)\leq\frac{1}{1+\Gamma(1+\gamma)^{-1}t}, \quad t>0.
\end{equation*}
Then it follows 
\begin{equation*}
    0< \mathbb{E}_{\gamma,1}(-t) <1, \quad t>0.
\end{equation*}
If $\gamma=1$, we know that $0<e^{-t}<1$, when $t>0$.

\begin{proof}
Solution of the Problem \ref{Problem 1} can be found by applying the Dunkl transform $\mathcal{F}_\alpha$ \eqref{Dunkl transform} and using formula \eqref{Commutation formula} to the equation \eqref{FracEq1} and the initial condition \eqref{FracCon1}. Thus, we have
\begin{equation}\label{FracEq2}
    \mathbb{D}_{0^+,t}^\gamma\widehat{u}(t,\lambda) + \frac{m+\lambda^2}{1 + a\lambda^2}\widehat{u}(t,\lambda) = \frac{\widehat{f}(t,\lambda)}{1 + a\lambda^2}, \quad \lambda\in\mathbb{R},
\end{equation}
and
\begin{equation}\label{FracCon2}
    \widehat{u}(0,\lambda)=\widehat{g}(\lambda), \quad \lambda\in\mathbb{R},
\end{equation}
where $\widehat{u}(\cdot,\lambda)$ is an unknown function. Let $0<\gamma\leq1$. The equation \eqref{FracEq2} is a ordinary differential equation respect to $t$, then by solving the equation \eqref{FracEq2} under the initial condition \eqref{FracCon2} (see \cite[p. 231, ex. 4.9]{KST06}), we obtain
\begin{multline} \label{DunSolution}
    \widehat{u}(t,\lambda) = \widehat{g}(\lambda)\mathbb{E}_{\gamma,1}\left(-\frac{m+\lambda^2}{1+a\lambda^2}t^\gamma\right)\\
    + \int_0^t(t-\tau)^{\gamma-1} \mathbb{E}_{\gamma,\gamma}\left(-\frac{m+\lambda^2}{1+a\lambda^2}(t-\tau)^\gamma\right)\frac{\widehat{f}(\tau,\lambda)}{1+a\lambda^2}d\tau,
\end{multline}
where $\mathbb{E}_{\gamma,1}(z)$ is the classical Mittag-Leffler function and $\mathbb{E}_{\gamma,\gamma}(z)$ is the Mittag-Leffler type function. Consequently, solution of the Problem \ref{Problem 1} is
\begin{equation*}
    u(t,x) = \int_\mathbb{R}\int_\mathbb{R} \mathbb{E}_{\gamma,1}\left(-\frac{m+\lambda^2}{1+a\lambda^2}t^\gamma\right)g(y)D_\alpha(ix\lambda)D_\alpha(-iy\lambda)d\mu_\alpha(y)d\mu_\alpha(\lambda)
\end{equation*}
\begin{multline*}
    + \int_\mathbb{R}\int_\mathbb{R}\int_0^t  (t-\tau)^{\gamma-1}\mathbb{E}_{\gamma,\gamma}\left(-\frac{m+\lambda^2}{1+a\lambda^2}(t-\tau)^\gamma\right)\frac{f(\tau,y)}{1+a\lambda^2}\\
    \times D_\alpha(ix\lambda)D_\alpha(-iy\lambda)d\tau d\mu_\alpha(y)d\mu_\alpha(\lambda),
\end{multline*}
here we have used the Fubini's theorem and the inverse Dunkl transform $\mathcal{F}_\alpha^{-1}$ \eqref{Invers Dunkl Transform} to \eqref{DunSolution}.  

By using the property
\begin{equation*}
    \frac{d}{d\tau}\left(\mathbb{E}_{\gamma,1}(c \tau^\gamma)\right)=c\tau^{\gamma-1}\mathbb{E}_{\gamma,\gamma}(c \tau^\gamma), \quad c =\text{constant},
\end{equation*}
of the Mittag-Leffler function, we obtain
\begin{equation*}
    \frac{d}{d\tau}\left(\mathbb{E}_{\gamma,1}\left(-\frac{m+\lambda^2}{1+a\lambda^2} (t-\tau)^\gamma\right)\right)=\frac{m+\lambda^2}{1+a\lambda^2} (t-\tau)^{\gamma-1}\mathbb{E}_{\gamma,\gamma}\left(-\frac{m+\lambda^2}{1+a\lambda^2} (t-\tau)^\gamma\right)
\end{equation*}
and 
\begin{multline} \label{Solutionf}
    \widehat{u}(t,\lambda) = \widehat{g}(\lambda)\mathbb{E}_{\gamma,1}\left(-\frac{m+\lambda^2}{1+a\lambda^2}t^\gamma\right)\\
    +\frac{1}{m+\lambda^2}\int_0^t \frac{d}{d\tau}\left(\mathbb{E}_{\gamma,1}\left(-\frac{m+\lambda^2}{1+a\lambda^2} (t-\tau)^\gamma\right)\right)\widehat{f}(\tau,\lambda)d\tau
\end{multline}
\begin{multline*}
    =\widehat{g}(\lambda)\mathbb{E}_{\gamma,1}\left(-\frac{m+\lambda^2}{1+a\lambda^2}t^\gamma\right)+\frac{\widehat{f}(t,\lambda)}{m+\lambda^2}-\frac{1}{m+\lambda^2}\mathbb{E}_{\gamma,1}\left(-\frac{m+\lambda^2}{1+a\lambda^2}t^\gamma\right)\widehat{f}(0,\lambda)\\
    -\frac{1}{m+\lambda^2}\int_0^t \mathbb{E}_{\gamma,1}\left(-\frac{m+\lambda^2}{1+a\lambda^2}(t-\tau)^\gamma\right)\frac{d}{d\tau}\widehat{f}(\tau,\lambda)d\tau
\end{multline*}
by using the Integration by Parts and $\mathbb{E}_{\gamma,1}(0)=1$.

Let $0<\gamma<1$. We assume that $g\in \mathcal{H}_\alpha^2(\mathbb{R},\mu_\alpha)$ and $f\in C^1([0,T],L^2(\mathbb{R},\mu_\alpha))$. Then for the function $u$ we have the
following estimate
\begin{equation*}
    \|u(t,\cdot)\|_{\mathcal{H}_\alpha^2}^2 = \int_\mathbb{R} (1+\lambda^2)^2|\widehat{u}(t,\lambda)|^2 d\mu_\alpha(\lambda)
\end{equation*}
\begin{align*}
    &\lesssim \int_\mathbb{R} (1+\lambda^2)^2\left| \widehat{g}(\lambda)\mathbb{E}_{\gamma,1}\left(-\frac{m+\lambda^2}{1+a\lambda^2}t^\gamma\right) \right|^2 d\mu_\alpha(\lambda)+\int_\mathbb{R} (1+\lambda^2)^2\left| \frac{\widehat{f}(t,\lambda)}{m+\lambda^2} \right|^2 d\mu_\alpha(\lambda)\\
    &+\int_\mathbb{R} (1+\lambda^2)^2\left| \frac{1}{m+\lambda^2}\mathbb{E}_{\gamma,1}\left(-\frac{m+\lambda^2}{1+a\lambda^2}t^\gamma\right)\widehat{f}(0,\lambda) \right|^2 d\mu_\alpha(\lambda)\\
    &+ \int_\mathbb{R} (1+\lambda^2)^2\left| \frac{1}{m+\lambda^2}\int_0^t \mathbb{E}_{\gamma,1}\left(-\frac{m+\lambda^2}{1+a\lambda^2}(t-\tau)^\gamma\right)\frac{d}{d\tau}\widehat{f}(\tau,\lambda)d\tau \right|^2 d\mu_\alpha(\lambda)\\
    &\lesssim \int_\mathbb{R} (1+\lambda^2)^2\left| \widehat{g}(\lambda) \right|^2 d\mu_\alpha(\lambda) + \int_\mathbb{R} \left| \widehat{f}(t,\lambda) \right|^2 d\mu_\alpha(\lambda)\\
    &+ \int_\mathbb{R} \left| \widehat{f}(0,\lambda) \right|^2 d\mu_\alpha(\lambda) + \int_\mathbb{R} \left( \int_0^t \left|\frac{d}{d\tau}\widehat{f}(\tau,\lambda)\right|d\tau \right)^2 d\mu_\alpha(\lambda)\\
    &\lesssim \|g\|_{\mathcal{H}_\alpha^2}^2 + \|f(t,\cdot)\|_{2, \alpha}^2 + \|f(0,\cdot)\|_{2, \alpha}^2 + \int_0^T \|\frac{d}{d t}f(t,\cdot)\|_{2, \alpha}^2 dt,
\end{align*}
where we have used the Fubini's theorem and $U\lesssim W$, which denotes $U\leq CW$ for some positive constant $C$ independent of $U$ and $W$. Thus, we obtain
\begin{equation*}
    \|u\|_{C([0,T],\mathcal{H}_\alpha^2(\mathbb{R},\mu_\alpha))}^2 \lesssim \|g\|_{\mathcal{H}_\alpha^2}^2 + \|f\|_{C^1([0,T],L^2(\mathbb{R},\mu_\alpha))}^2 < +\infty.
\end{equation*}

Now, for $\mathcal{D}_{0^+,t}^\gamma u$ we have 
\begin{align*}
    \|\mathcal{D}_{0^+,t}^\gamma u(t,\cdot)\|_{2,\alpha}^2 &= \|\mathcal{F}_{\alpha}\left(\mathcal{D}_{0^+,t}^\gamma u(t,\cdot)\right)\|_{2,\alpha}^2 = \|\mathcal{D}_{0^+,t}^\gamma \widehat{u}(t,\cdot)\|_{2,\alpha}^2\\
    &= \int_\mathbb{R} \left| \frac{\widehat{f}(t,\lambda)}{1 + a\lambda^2} - \frac{m+\lambda^2}{1 + a\lambda^2}\widehat{u}(t,\lambda) \right|^2 d\mu_\alpha(\lambda)\\
    &\lesssim \int_\mathbb{R}\left| \widehat{f}(t,\lambda) \right|^2 d\mu_\alpha(\lambda) + \int_\mathbb{R}(1+\lambda^2)^2\left| \widehat{u}(t,\lambda)\right|^2 d\mu_\alpha(\lambda)\\
    &= \|f(t,\cdot)\|_{2,\alpha}^2 + \|u(t,\cdot)\|_{\mathcal{H}_\alpha^2}^2
\end{align*}
by using \eqref{PTh}. Consequently, it gives us
\begin{equation*}
    \|\mathcal{D}_{0^+,t}^\gamma u\|_{C([0,T],L^2(\mathbb{R},\mu_\alpha))}^2 \lesssim \|f\|_{C([0,T],L^2(\mathbb{R},\mu_\alpha))}^2 + \|u\|_{C([0,T],\mathcal{H}_\alpha^2(\mathbb{R},\mu_\alpha))}^2 < +\infty.
\end{equation*}
Finally, using Definition \ref{Defn} we obtain $u\in C^\gamma([0,T], L^2(\mathbb{R},\mu_\alpha))$. 

Let $\gamma=1$. We assume that $g\in \mathcal{H}_\alpha^2(\mathbb{R},\mu_{\alpha})$ and $f\in C([0,T],L^2(\mathbb{R},\mu_{\alpha}))$. Then let us estimate the function $u$ as follows 
\begin{equation*}
    \|u(t,\cdot)\|_{\mathcal{H}_\alpha^2}^2 = \int_\mathbb{R}|(1+\lambda^2)\widehat{u}(t,\lambda)|^2d\mu_\alpha(\lambda)
\end{equation*}
\begin{align*}
    &= \int_\mathbb{R}\left|\int_0^t \frac{1+\lambda^2}{1+a\lambda^2} \widehat{f}(\tau,\lambda)e^{-\frac{m+\lambda^2}{1 + a\lambda^2}(t-\tau)}d\tau + (1+\lambda^2)\widehat{g}(\lambda)e^{-\frac{m+\lambda^2}{1 + a\lambda^2}t}\right|^2d\mu_\alpha(\lambda)\\
    &\lesssim \int_\mathbb{R}\left|\frac{1+\lambda^2}{1+a\lambda^2}\int_0^t \widehat{f}(\tau,\lambda)e^{-\frac{m+\lambda^2}{1 + a\lambda^2}(t-\tau)}d\tau\right|^2d\mu_\alpha(\lambda) + \int_\mathbb{R}\left|(1+\lambda^2)\widehat{g}(\lambda)e^{-\frac{m+\lambda^2}{1 + a\lambda^2}t}\right|^2d\mu_\alpha(\lambda)\\
    &\lesssim \int_\mathbb{R}\left(\int_0^T |\widehat{f}(t,\lambda)|dt\right)^2d\mu_\alpha(\lambda) + \int_\mathbb{R}|(1+\lambda^2)\widehat{g}(\lambda)|^2d\mu_\alpha(\lambda)\\
    &\lesssim \int_0^T\|\widehat{f}(t,\cdot)\|_{2,\alpha}^2dt + \|g\|_{\mathcal{H}_\alpha^2}^2\\ 
    &= \int_0^T\|f(t,\cdot)\|_{2,\alpha}^2dt + \|g\|_{\mathcal{H}_\alpha^2}^2
\end{align*}
by using H\"{o}lder's inequality, Fubini's theorem and \eqref{PTh}. Then
\begin{equation*}
    \|u\|_{C([0,T],\mathcal{H}_\alpha^2(\mathbb{R},\mu_{\alpha}))}^2 \lesssim \|f\|_{C([0,T],L^2(\mathbb{R},\mu_{\alpha}))}^2 + \|g\|_{\mathcal{H}_\alpha^2}^2 < +\infty.
\end{equation*}
Let us estimate $u_t$ as follows
\begin{multline*}
    \|u_t(t,\cdot)\|_{2,\alpha}^2 = \|\widehat{u}_t(t,\cdot)\|_{2,\alpha}^2 = \int_\mathbb{R}|\widehat{u}_t(t,\lambda)|^2d\mu_{\alpha}(\lambda)\\
    = \int_\mathbb{R} \left| \frac{\widehat{f}(t,\lambda)}{1 + a\lambda^2} - \frac{m+\lambda^2}{1 + a\lambda^2}\widehat{u}(t,\lambda) \right|^2 d\mu_\alpha(\lambda)\\
    \lesssim \|f(t,\cdot)\|_{2,\alpha}^2 + \|u(t,\cdot)\|_{\mathcal{H}_\alpha^2}^2.
\end{multline*}
Thus 
\begin{align*}
    \|u_t\|_{C([0,T],L^2(\mathbb{R},\mu_\alpha))}^2 &\lesssim \|f\|_{C([0,T],L^2(\mathbb{R},\mu_\alpha))}^2 + \|u\|_{C([0,T],\mathcal{H}_\alpha^2(\mathbb{R},\mu_\alpha))}^2 < +\infty.
\end{align*}
The existence is proved.

Now, we are going to prove uniqueness of the solution. Suppose that there are two solutions $u_1$ and $u_2$ of the Problem \ref{Problem 1}. Denote
\begin{equation*}
    u(t,x) = u_1(t,x) - u_2(t,x).
\end{equation*}
Then the function $u$ is a solution of the problem
\begin{equation*}
    \begin{cases}
    \mathbb{D}_{0^+,t}^\gamma \left(u(t,x) - a\Lambda_{\alpha,x}^2 u(t,x)\right) - \Lambda_{\alpha,x}^2 u(t,x) + m u(t,x) = 0,\\
    u(0,\lambda)=0.
    \end{cases}
\end{equation*}
Then by applying the Dunkl transform $\mathcal{F}_\alpha$ \eqref{Dunkl transform}, we obtain
\begin{equation*}
    \begin{cases}
    \mathbb{D}_{0^+,t}^\gamma\widehat{u}(t,\lambda) + \frac{m+\lambda^2}{1 + a\lambda^2}\widehat{u}(t,\lambda) = 0,\\
    \widehat{u}(0,\lambda)=0.
    \end{cases}
\end{equation*}
Above equation has a trivial solution (see \cite[p. 231, ex. 4.9]{KST06}), i.e. $\widehat{u}(t,\lambda)\equiv0$. Hence, uniqueness of the solution is proved.
\end{proof}

\section{Inverse Source Problem}
\label{Sect.In}

\subsection{Problem} 
Here we deal with the following ISP for the equation \eqref{EQ:Heat Dunkl}.

\begin{problem} \label{Problem 2}
Let $0<\gamma\leq1$. Our aim to find the couple of functions $(u,f)$ satisfying the equation 
\begin{equation}
\label{FracInEq1}
    \mathbb{D}_{0^+,t}^\gamma \left(u(t,x) - a\Lambda_{\alpha,x}^2 u(t,x)\right) - \Lambda_{\alpha,x}^2 u(t,x) + m u(t,x) = f(x), \quad (t,x)\in Q_T,
\end{equation}
under the conditions
\begin{equation}
\label{FracInCon1}
    u(0,x)=\phi(x), \quad x\in\mathbb{R}
\end{equation}
and
\begin{equation}
\label{FracInCon2}
    u(T,x)=\psi(x), \quad x\in\mathbb{R},
\end{equation}
where $\phi$ and $\psi$ is sufficiently smooth functions, $\Lambda_\alpha$ is the Dunkl operator.
\end{problem}

We assume that $0<\gamma\leq1$. Then a generalised  solution of Problem \ref{Problem 2} is a pair of functions $(u,f)$ satisfying the above problem such that $u\in C^\gamma([0,T],L^2(\mathbb{R},\mu_\alpha))\cap C([0,T],\mathcal{H}_\alpha^2(\mathbb{R},\mu_\alpha))$ and $f\in L^2(\mathbb{R},\mu_\alpha)$.

\begin{theorem} \label{FracInverseTheorem}
Let $0<\gamma\leq1$. We assume that $\psi,\phi\in \mathcal{H}_\alpha^2(\mathbb{R},\mu_\alpha)$. Then the solution of Problem \ref{Problem 2} exists, is unique, and can be written by the expressions
\begin{multline*}
    f(x) = \int_\mathbb{R}\int_\mathbb{R} (m+\lambda^2)\frac{\psi(y) - \phi(y)\mathbb{E}_{\gamma,1}\left(-\frac{m+\lambda^2}{1 + a\lambda^2}T^\gamma\right)}{1 - \mathbb{E}_{\gamma,1}\left(-\frac{m+\lambda^2}{1 + a\lambda^2}T^\gamma\right)}\\
    \times D_\alpha(ix\lambda)D_\alpha(-iy\lambda)d\mu_\alpha(y)d\mu_\alpha(\lambda)
\end{multline*}
and 
\begin{multline*}
    u(t,x) =  \int_\mathbb{R}\int_\mathbb{R}\frac{1 - \mathbb{E}_{\gamma,1}\left(-\frac{m+\lambda^2}{1 + a\lambda^2}t^\gamma\right)}{1 - \mathbb{E}_{\gamma,1}\left(-\frac{m+\lambda^2}{1 + a\lambda^2}T^\gamma\right)}\psi(y) D_\alpha(ix\lambda)D_\alpha(-iy\lambda)d\mu_\alpha(y)d\mu_\alpha(\lambda)\\
    -\int_\mathbb{R}\int_\mathbb{R}\frac{\mathbb{E}_{\gamma,1}\left(-\frac{m+\lambda^2}{1 + a\lambda^2}T^\gamma\right)-\mathbb{E}_{\gamma,1}\left(-\frac{m+\lambda^2}{1 + a\lambda^2}t^\gamma\right)}{1 - \mathbb{E}_{\gamma,1}\left(-\frac{m+\lambda^2}{1 + a\lambda^2}T^\gamma\right)}\phi(y) D_\alpha(ix\lambda)D_\alpha(-iy\lambda)d\mu_\alpha(y)d\mu_\alpha(\lambda).
\end{multline*}
\end{theorem}

\begin{proof}
Repeating the method applied in the previous section, we aim at finding a solution to Problem \ref{Problem 2} by applying the Dunkl transform $\mathcal{F}_\alpha$ \eqref{Dunkl transform} to the equation \eqref{FracInEq1} and the conditions \eqref{FracInCon1} and \eqref{FracInCon2}. It gives us
\begin{equation}
\label{FracInEq2}
    \mathbb{D}_{0^+,t}^\gamma\widehat{u}(t,\lambda) + \frac{m+\lambda^2}{1 + a\lambda^2}\widehat{u}(t,\lambda) = \frac{\widehat{f}(\lambda)}{1 + a\lambda^2}, \quad \lambda\in\mathbb{R},
\end{equation}
and
\begin{equation}
\label{FracInCon3}
    \widehat{u}(0,\lambda)=\widehat{\phi}(\lambda), \quad \lambda\in\mathbb{R},
\end{equation}
\begin{equation}
\label{FracInCon4}
    \widehat{u}(T,\lambda)=\widehat{\psi}(\lambda), \quad \lambda\in\mathbb{R},
\end{equation}
where $\widehat{u}(t,\lambda)$ and $\widehat{f}(\lambda)$ are unknown. 

Let $0<\gamma\leq1$. Using expression \eqref{Solutionf} we can find solution of the equation \eqref{FracInEq2} with initial condition \eqref{FracInCon3}, given by
\begin{equation}
\label{GeneralInSolution}
    \widehat{u}(t,\lambda) = \frac{\widehat{f}(\lambda)}{m+\lambda^2} + \left(\widehat{\phi}(\lambda) - \frac{\widehat{f}(\lambda)}{m+\lambda^2}\right)\mathbb{E}_{\gamma,1}\left(-\frac{m+\lambda^2}{1 + a\lambda^2}t^\gamma\right),
\end{equation}
where $\widehat{f}(\lambda)$ is unknown and $\mathbb{E}_{\gamma,1}\left(z\right)$ is the Mittag-Leffler function. Then applying the condition \eqref{FracInCon4} to the expression \eqref{GeneralInSolution}, one obtains
\begin{equation*}
    \widehat{u}(T,\lambda) = \frac{\widehat{f}(\lambda)}{m+\lambda^2} + \left(\widehat{\phi}(\lambda) - \frac{\widehat{f}(\lambda)}{m+\lambda^2}\right)\mathbb{E}_{\gamma,1}\left(-\frac{m+\lambda^2}{1 + a\lambda^2}T^\gamma\right) = \widehat{\psi}(\lambda)
\end{equation*}
Thus, we can find unknown $\widehat{f}(\lambda)$ as following 
\begin{equation} \label{Dunklsolution1}
    \widehat{f}(\lambda) = (m+\lambda^2)\frac{\widehat{\psi}(\lambda) - \widehat{\phi}(\lambda)\mathbb{E}_{\gamma,1}\left(-\frac{m+\lambda^2}{1 + a\lambda^2}T^\gamma\right)}{1 - \mathbb{E}_{\gamma,1}\left(-\frac{m+\lambda^2}{1 + a\lambda^2}T^\gamma\right)}.
\end{equation}
Consequently, by substituting $\widehat{f}(\lambda)$ into \eqref{GeneralInSolution}, we have
\begin{equation}
\label{Dunklsolution2}
    \widehat{u}(t,\lambda)= \frac{1 - \mathbb{E}_{\gamma,1}\left(-\frac{m+\lambda^2}{1 + a\lambda^2}t^\gamma\right)}{1 - \mathbb{E}_{\gamma,1}\left(-\frac{m+\lambda^2}{1 + a\lambda^2}T^\gamma\right)}\widehat{\psi}(\lambda) - \frac{\mathbb{E}_{\gamma,1}\left(-\frac{m+\lambda^2}{1 + a\lambda^2}T^\gamma\right) - \mathbb{E}_{\gamma,1}\left(-\frac{m+\lambda^2}{1 + a\lambda^2}t^\gamma\right)}{1 - \mathbb{E}_{\gamma,1}\left(-\frac{m+\lambda^2}{1 + a\lambda^2}T^\gamma\right)}\widehat{\phi}(\lambda).
\end{equation}
Then, one obtains the solution to Problem \ref{Problem 2}, which is a pair of functions $(u,f)$ given by the formulas
\begin{multline*}
    f(x) = \int_\mathbb{R}\int_\mathbb{R} (m+\lambda^2)\frac{\psi(y) - \phi(y)\mathbb{E}_{\gamma,1}\left(-\frac{m+\lambda^2}{1 + a\lambda^2}T^\gamma\right)}{1 - \mathbb{E}_{\gamma,1}\left(-\frac{m+\lambda^2}{1 + a\lambda^2}T^\gamma\right)}\\
    \times D_\alpha(ix\lambda)D_\alpha(-iy\lambda)d\mu_\alpha(y)d\mu_\alpha(\lambda),
\end{multline*}
and 
\begin{multline*}
    u(t,x) =  \int_\mathbb{R}\int_\mathbb{R}\frac{1 - \mathbb{E}_{\gamma,1}\left(-\frac{m+\lambda^2}{1 + a\lambda^2}t^\gamma\right)}{1 - \mathbb{E}_{\gamma,1}\left(-\frac{m+\lambda^2}{1 + a\lambda^2}T^\gamma\right)}\psi(y) D_\alpha(ix\lambda)D_\alpha(-iy\lambda)d\mu_\alpha(y)d\mu_\alpha(\lambda)\\
    -\int_\mathbb{R}\int_\mathbb{R}\frac{\mathbb{E}_{\gamma,1}\left(-\frac{m+\lambda^2}{1 + a\lambda^2}T^\gamma\right)-\mathbb{E}_{\gamma,1}\left(-\frac{m+\lambda^2}{1 + a\lambda^2}t^\gamma\right)}{1 - \mathbb{E}_{\gamma,1}\left(-\frac{m+\lambda^2}{1 + a\lambda^2}T^\gamma\right)}\phi(y) D_\alpha(ix\lambda)D_\alpha(-iy\lambda)d\mu_\alpha(y)d\mu_\alpha(\lambda),
\end{multline*}
using the inverse Dunkl transform $\mathcal{F}_\alpha^{-1}$ \eqref{Invers Dunkl Transform} for the equations \eqref{Dunklsolution1} and \eqref{Dunklsolution2}, respectively.

Let $\psi,\phi\in \mathcal{H}_\alpha^2(\mathbb{R},\mu_\alpha)$. Then, we have the following estimate
\begin{multline*}
    \|f\|_{2,\alpha}^2 = \|\widehat{f}\|_{2,\alpha}^2 = \int_\mathbb{R} |\widehat{f}(\lambda)|^2 d\mu_\alpha(\lambda)\\
    = \int_\mathbb{R}(m+\lambda^2)^2 \left|\frac{\widehat{\psi}(\lambda)-\widehat{\phi}(\lambda)\mathbb{E}_{\gamma,1}\left(-\frac{m+\lambda^2}{1 + a\lambda^2}T^\gamma\right)}{1 - \mathbb{E}_{\gamma,1}\left(-\frac{m+\lambda^2}{1 + a\lambda^2}T^\gamma\right)} \right|^2 d\mu_\alpha(\lambda)
\end{multline*}
\begin{multline*}
    \lesssim \int_\mathbb{R}\frac{(m+\lambda^2)^2}{\left(1 - \mathbb{E}_{\gamma,1}\left(-\frac{m+\lambda^2}{1 + a\lambda^2}T^\gamma\right)\right)^2}|\widehat{\psi}(\lambda)|^2d\mu_\alpha(\lambda)\\
    + \int_\mathbb{R}(m+\lambda^2)^2\left(\frac{\mathbb{E}_{\gamma,1}\left(-\frac{m+\lambda^2}{1 + a\lambda^2}T^\gamma\right)}{1 - \mathbb{E}_{\gamma,1}\left(-\frac{m+\lambda^2}{1 + a\lambda^2}T^\gamma\right)}\right)^2|\widehat{\psi}(\lambda)|^2d\mu_\alpha(\lambda)\\
    \lesssim \int_\mathbb{R}(1+\lambda^2)^2|\widehat{\psi}(\lambda)|^2d\mu_\alpha(\lambda) + \int_\mathbb{R}(1+\lambda^2)^2|\widehat{\phi}(\lambda)|^2d\mu_\alpha(\lambda),
\end{multline*}
by using \eqref{PTh} and properties of Mittag-Leffler functions. Thus, we have
\begin{equation*}
    \|f\|_{2,\alpha}^2 \lesssim \|\psi\|_{\mathcal{H}_\alpha^2}^2 + \|\phi\|_{\mathcal{H}_\alpha^2}^2 < +\infty.
\end{equation*}
For $u(t,x)$, one arrives 
\begin{align*}
    \|u(t,\cdot)\|_{\mathcal{H}_\alpha^2}^2 &= \int_\mathbb{R} (1+\lambda^2)^2|\widehat{u}(t,\lambda)|^2 d\mu_\alpha(\lambda)\\
    &\lesssim \int_\mathbb{R} (1+\lambda^2)^2\left| \frac{1 - \mathbb{E}_{\gamma,1}\left(-\frac{m+\lambda^2}{1 + a\lambda^2}t^\gamma\right)}{1 - \mathbb{E}_{\gamma,1}\left(-\frac{m+\lambda^2}{1 + a\lambda^2}T^\gamma\right)}\widehat{\psi}(\lambda)\right|^2 d\mu_\alpha(\lambda)\\
    &+ \int_\mathbb{R} (1+\lambda^2)^2\left| \frac{\mathbb{E}_{\gamma,1}\left(-\frac{m+\lambda^2}{1 + a\lambda^2}T^\gamma\right) - \mathbb{E}_{\gamma,1}\left(-\frac{m+\lambda^2}{1 + a\lambda^2}t^\gamma\right)}{1 - \mathbb{E}_{\gamma,1}\left(-\frac{m+\lambda^2}{1 + a\lambda^2}T^\gamma\right)}\widehat{\phi}(\lambda) \right|^2 d\mu_\alpha(\lambda)\\
    &\lesssim \|\psi\|_{\mathcal{H}_\alpha^2}^2 + \|\phi\|_{\mathcal{H}_\alpha^2}^2.
\end{align*}
Consequently, it gives 
\begin{equation*}
    \|u\|_{C([0,T],\mathcal{H}_\alpha^2(\mathbb{R},\mu_\alpha))}^2 \lesssim \|\psi\|_{\mathcal{H}_\alpha^2}^2 + \|\phi\|_{\mathcal{H}_\alpha^2}^2 < +\infty.
\end{equation*}
For $u_t(t,\cdot)$, we have
\begin{equation*}
    \mathbb{D}_{0^+,t}^\gamma\widehat{u}(t,\lambda) = \frac{\widehat{f}(\lambda)}{1 + a\lambda^2} - \frac{m+\lambda^2}{1 + a\lambda^2}\widehat{u}(t,\lambda),
\end{equation*}
rewriting the equation \eqref{FracInEq2}. Then by using \eqref{PTh}, one obtains
\begin{multline*}
    \|\mathbb{D}_{0^+,t}^\gamma u(t,\cdot)\|_{2,\alpha}^2 = \|\mathcal{F}_{\alpha}\left(\mathbb{D}_{0^+,t}^\gamma u(t,\cdot)\right)\|_{2,\alpha}^2 \\
    = \|\mathbb{D}_{0^+,t}^\gamma \widehat{u}(t,\cdot)\|_{2,\alpha}^2 = \int_\mathbb{R} \left|\frac{\widehat{f}(\lambda)}{1 + a\lambda^2} - \frac{m+\lambda^2}{1 + a\lambda^2}\widehat{u}(t,\lambda)\right|^2 d\mu_\alpha(\lambda)\\
    \lesssim \|f\|_{2,\alpha}^2 + \|u(t,\cdot)\|_{\mathcal{H}_\alpha^2}^2.
\end{multline*}
Thus,
\begin{equation*}
    \|\mathbb{D}_{0^+,t}^\gamma u\|_{C([0,T],L^2(\mathbb{R},\mu_\alpha))}^2 \lesssim \|f\|_{2,\alpha}^2 + \|u\|_{C([0,T],\mathcal{H}_\alpha^2(\mathbb{R},\mu_\alpha))}^2<+\infty.
\end{equation*}
The existence is proved.

Now, we are going to prove the uniqueness of the solution. Suppose that there are two solutions $(u_1,f_1)$ and $(u_2,f_2)$ of Problem \ref{Problem 2}. Denote
\begin{equation*}
    u(t,x) = u_1(t,x) - u_2(t,x)
\end{equation*}
and
\begin{equation*}
    f(x) = f_1(x) - f_2(x).
\end{equation*}
Then the functions $u$ and $f$ satisfy 
\begin{equation*}
    \begin{cases}
    \mathbb{D}_{0^+,t}^\gamma \left(u(t,x) - a\Lambda_{\alpha,x}^2 u(t,x)\right) - \Lambda_{\alpha,x}^2 u(t,x) + m u(t,x) = f(x),\\
    u(0,\lambda)=0,\\
    u(T,\lambda)=0.
    \end{cases}
\end{equation*}
Then by applying the Dunkl transform $\mathcal{F}_\alpha$ \eqref{Dunkl transform}, one has
\begin{equation*}
    \begin{cases}
    \mathbb{D}_{0^+,t}^\gamma\widehat{u}(t,\lambda) + \frac{m + \lambda^2}{1 + a\lambda^2}\widehat{u}(t,\lambda) = \frac{\widehat{f}(\lambda)}{1 + a\lambda^2},\\
    \widehat{u}(0,\lambda)=0,\\
    \widehat{u}(T,\lambda)=0.
    \end{cases}
\end{equation*}
Via our calculation above, we can see that the problem has a trivial solution, i.e. $\widehat{u}(t,\lambda)\equiv0$ and $\widehat{f}(\lambda)\equiv0$. Hence, the uniqueness of the solution is proved.
\end{proof}

\subsection{Stability analysis of ISP}

In this subsection we study stability of the solution $(u,f)$ which is given by the formulas \eqref{Dunklsolution1} and \eqref{Dunklsolution2}, of Problem \ref{Problem 2}.

\begin{theorem} Let $(u,f)$ and $(u_d,f_d)$ be solutions to Problem \ref{Problem 2} corresponding to the
data $(\phi,\psi)$ and its small perturbation $(\phi_d,\psi_d)$, respectively. Then the solution of Problem \ref{Problem 2} depends continuously on these data, namely, we have
\begin{equation*}
    \|u-u_d\|_{C([0,T],\mathcal{H}_\alpha^2(\mathbb{R},\mu_\alpha))}^2 \lesssim \|\psi-\psi_d\|_{\mathcal{H}_\alpha^2}^2 + \|\phi-\phi_d\|_{\mathcal{H}_\alpha^2}^2
\end{equation*}
and
\begin{equation*}
    \|f-f_d\|_{2,\alpha}^2 \lesssim \|\psi-\psi_d\|_{\mathcal{H}_\alpha^2}^2 + \|\phi-\phi_d\|_{\mathcal{H}_\alpha^2}^2.
\end{equation*}
\end{theorem}

\begin{proof}
From the definition of the Dunkl transform
\begin{equation*}
    \mathcal{F}_{\alpha}(u(t,\cdot))(\lambda) = \widehat{u}(t,\lambda) = \int_\mathbb{R} u(t,x)D_{\alpha}(-i x\lambda) d\mu_\alpha(x),
\end{equation*}
we have
\begin{align*}
    \mathcal{F}_{\alpha}(u(t,\cdot)-u_d(t,\cdot))(\lambda) &= \int_\mathbb{R} (u(t,x)-u_d(t,x))D_{\alpha}(-i x\lambda) d\mu_\alpha(x)\\
    &= \int_\mathbb{R} u(t,x)D_{\alpha}(-i x\lambda) d\mu_\alpha(x) - \int_\mathbb{R} u_d(t,x)D_{\alpha}(-i x\lambda) d\mu_\alpha(x)\\
    &= \mathcal{F}_{\alpha}(u(t,\cdot))(\lambda) - \mathcal{F}_{\alpha}(u_d(t,\cdot))(\lambda)\\
    &=\widehat{u}(t,\lambda)-\widehat{u}_d(t,\lambda),
\end{align*}
using the property of the integral.

\begin{equation*}
    \|u(t,\cdot)-u_d(t,\cdot)\|_{\mathcal{H}_\alpha^2}^2 = \int_\mathbb{R} (1+\lambda^2)^2|\widehat{u}(t,\lambda)-\widehat{u}_d(t,\lambda)|^2 d\mu_\alpha(\lambda)
\end{equation*}
\begin{multline*}
    =\int_\mathbb{R} (1+\lambda^2)^2\Bigg|\frac{1 - \mathbb{E}_{\gamma,1}\left(-\frac{m+\lambda^2}{1 + a\lambda^2}t^\gamma\right)}{1 - \mathbb{E}_{\gamma,1}\left(-\frac{m+\lambda^2}{1 + a\lambda^2}T^\gamma\right)}\left(\widehat{\psi}(\lambda)-\widehat{\psi}_d(\lambda)\right)\\
    -\frac{\mathbb{E}_{\gamma,1}\left(-\frac{m+\lambda^2}{1 + a\lambda^2}T^\gamma\right) - \mathbb{E}_{\gamma,1}\left(-\frac{m+\lambda^2}{1 + a\lambda^2}t^\gamma\right)}{1 - \mathbb{E}_{\gamma,1}\left(-\frac{m+\lambda^2}{1 + a\lambda^2}T^\gamma\right)}\left(\widehat{\phi}(\lambda)-\widehat{\phi}_d(\lambda)\right)\bigg|^2d\mu_\alpha(\lambda)
\end{multline*}
\begin{align*}
    &\lesssim \int_\mathbb{R} (1+\lambda^2)^2\left|\frac{1 - \mathbb{E}_{\gamma,1}\left(-\frac{m+\lambda^2}{1 + a\lambda^2}t^\gamma\right)}{1 - \mathbb{E}_{\gamma,1}\left(-\frac{m+\lambda^2}{1 + a\lambda^2}T^\gamma\right)}\left(\widehat{\psi}(\lambda)-\widehat{\psi}_d(\lambda)\right)\right|^2d\mu_\alpha(\lambda)\\
    &+ \int_\mathbb{R} (1+\lambda^2)^2\left|\frac{\mathbb{E}_{\gamma,1}\left(-\frac{m+\lambda^2}{1 + a\lambda^2}T^\gamma\right) - \mathbb{E}_{\gamma,1}\left(-\frac{m+\lambda^2}{1 + a\lambda^2}t^\gamma\right)}{1 - \mathbb{E}_{\gamma,1}\left(-\frac{m+\lambda^2}{1 + a\lambda^2}T^\gamma\right)}\left(\widehat{\phi}(\lambda)-\widehat{\phi}_d(\lambda)\right)\right|^2d\mu_\alpha(\lambda)\\
    &\lesssim \int_\mathbb{R} (1+\lambda^2)^2\left|\widehat{\psi}(\lambda)-\widehat{\psi}_d(\lambda)\right|^2d\mu_\alpha(\lambda) + \int_\mathbb{R} (1+\lambda^2)^2\left|\widehat{\phi}(\lambda)-\widehat{\phi}_d(\lambda)\right|^2d\mu_\alpha(\lambda)\\
    &\lesssim \|\psi-\psi_d\|_{\mathcal{H}_\alpha^2}^2 + \|\phi-\phi_d\|_{\mathcal{H}_\alpha^2}^2.
\end{align*}
Consequently,
\begin{equation*}
    \|u(t,\cdot)-u_d(t,\cdot)\|_{\mathcal{H}_\alpha^2}^2 \lesssim \|\psi-\psi_d\|_{\mathcal{H}_\alpha^2}^2 + \|\phi-\phi_d\|_{\mathcal{H}_\alpha^2}^2,
\end{equation*}
or
\begin{equation*}
    \|u-u_d\|_{C([0,T],\mathcal{H}_\alpha^2(\mathbb{R},\mu_\alpha))}^2 \lesssim \|\psi-\psi_d\|_{\mathcal{H}_\alpha^2}^2 + \|\phi-\phi_d\|_{\mathcal{H}_\alpha^2}^2.
\end{equation*}

By writing \eqref{Dunklsolution1} in the form
\begin{equation*}
    \widehat{f}(\lambda) = \frac{m+\lambda^2}{1 - \mathbb{E}_{\gamma,1}\left(-\frac{m+\lambda^2}{1 + a\lambda^2}T^\gamma\right)}\widehat{\psi}(\lambda) -\frac{(m+\lambda^2)\mathbb{E}_{\gamma,1}\left(-\frac{m+\lambda^2}{1 + a\lambda^2}T^\gamma\right)}{1 - \mathbb{E}_{\gamma,1}\left(-\frac{m}{1 + \lambda^2}T^\gamma\right)}\widehat{\phi}(\lambda),
\end{equation*}
one concludes
\begin{multline*}
    \|f-f_d\|_{2,\alpha}^2 = \|\widehat{f}-\widehat{f_d}\|_{2,\alpha}^2 = \int_\mathbb{R} |\widehat{f}(\lambda)-\widehat{f_d}(\lambda)|^2 d\mu_\alpha(\lambda)\\
    =\int_\mathbb{R} \bigg|\frac{m+\lambda^2}{1 - \mathbb{E}_{\gamma,1}\left(-\frac{m+\lambda^2}{1 + a\lambda^2}T^\gamma\right)}\left(\widehat{\psi}(\lambda)-\widehat{\psi}_d(\lambda)\right)\\
    - \frac{(m+\lambda^2)\mathbb{E}_{\gamma,1}\left(-\frac{m+\lambda^2}{1 + a\lambda^2}T^\gamma\right)}{1 - \mathbb{E}_{\gamma,1}\left(-\frac{m+\lambda^2}{1 + a\lambda^2}T^\gamma\right)}\left(\widehat{\phi}(\lambda)-\widehat{\phi}_d(\lambda)\right)\bigg|^2d\mu_\alpha(\lambda)
\end{multline*}
\begin{multline*}
    \lesssim \int_\mathbb{R} \left|\frac{m+\lambda^2}{1 - \mathbb{E}_{\gamma,1}\left(-\frac{m+\lambda^2}{1 + a\lambda^2}T^\gamma\right)}\left(\widehat{\psi}(\lambda)-\widehat{\psi}_d(\lambda)\right)\right|^2d\mu_\alpha(\lambda)\\
    + \int_\mathbb{R} \left|\frac{(m+\lambda^2)\mathbb{E}_{\gamma,1}\left(-\frac{m+\lambda^2}{1 + a\lambda^2}T^\gamma\right)}{1 - \mathbb{E}_{\gamma,1}\left(-\frac{m+\lambda^2}{1 + a\lambda^2}T^\gamma\right)}\left(\widehat{\phi}(\lambda)-\widehat{\phi}_d(\lambda)\right)\right|^2d\mu_\alpha(\lambda)\\
    \lesssim \|\psi-\psi_d\|_{\mathcal{H}_\alpha^2}^2 + \|\phi-\phi_d\|_{\mathcal{H}_\alpha^2}^2.
\end{multline*}
Thus,
\begin{equation*}
    \|f-f_d\|_{2,\alpha}^2 \lesssim \|\psi-\psi_d\|_{\mathcal{H}_\alpha^2}^2 + \|\phi-\phi_d\|_{\mathcal{H}_\alpha^2}^2.
\end{equation*}
We complete the proof.
\end{proof}

\subsection{Stability Test} Here we test one sample case for the subject of the stability of the solution pair. Let us consider the following inverse source problem for the pseudo-parabolic equation 
\begin{equation*}
    \frac{\partial}{\partial t}\left(u(t,x)-\frac{\partial^2}{\partial x^2}u(t,x)\right)-\frac{\partial^2}{\partial x^2}u(t,x)+u(t,x) = f(x), \quad 0<t<1, \quad x\in\mathbb{R}
\end{equation*}
with Dirichlet boundary conditions
\begin{equation*}
    u(0,x)=u(1,x)=0,
\end{equation*}
where $T=m=a=\gamma=1$, $\alpha=-\frac{1}{2}$, and $\phi(x)=\psi(x)=0$. 

By applying Theorem \ref{FracInverseTheorem}, with $\gamma=1$ and $\alpha=-\frac{1}{2}$: our operator in time is $\frac{d}{dt}$ and in space $\Lambda_{-\frac{1}{2},x}^2=\frac{d^2}{dx^2}$, we obtain the trivial solution pair:
\begin{equation*}
    u(t,x)\equiv0 \quad \text{and} \quad f(x)\equiv0.
\end{equation*}

Now we consider a perturbation of the previous problem in the following form
\begin{equation*}
    \frac{\partial}{\partial t}\left(u_d(t,x)-\frac{\partial^2}{\partial x^2}u_d(t,x)\right)-\frac{\partial^2}{\partial x^2}u_d(t,x)+u_d(t,x) = f_d(x), \quad 0<t<1, \quad x\in\mathbb{R}
\end{equation*}
 with conditions
\begin{equation*}
    u_d(0,x)=0, \quad \text{and} \quad u_d(1,x)=\epsilon\cdot\exp(-x^2), \quad \epsilon>0, \quad x\in\mathbb{R},
\end{equation*}
where $\phi_d(x)=0$, $\psi_d(x)=\epsilon\cdot\exp(-x^2)$ ($\psi_d\in\mathcal{H}_\alpha^2(\mathbb{R},\mu_\alpha)$). Then using Theorem \ref{FracInverseTheorem}, one obtains solution of the perturbation problem, expressed by
\begin{equation}\label{int1}
    u_d(t,x)=\frac{\epsilon}{2\sqrt{\pi}}\int_{-\infty}^\infty \frac{1-\exp(-t)}{1-\exp(-1)}\exp\left(-\frac{\lambda^2}{4}\right)\exp(ix\lambda)d\lambda
\end{equation}
and
\begin{equation}\label{int2}
    f_d(x)=\frac{\epsilon}{2\sqrt{\pi}}\int_{-\infty}^\infty \frac{1+\lambda^2}{1-\exp(-1)}\exp\left(-\frac{\lambda^2}{4}\right)\exp(ix\lambda)d\lambda.
\end{equation}
Integrals \eqref{int1} and \eqref{int2} are converges absolutely, because 
\begin{multline*}
    \left|u_d(t,x)\right|\leq\frac{\epsilon}{2\sqrt{\pi}}\int_{-\infty}^\infty\frac{1-\exp(-t)}{1-\exp(-1)}\exp\left(-\frac{\lambda^2}{4}\right)d\lambda\\
    =\frac{\epsilon}{2\sqrt{\pi}}\frac{1-\exp(-t)}{1-\exp(-1)}\int_{-\infty}^\infty\exp\left(-\frac{\lambda^2}{4}\right)d\lambda\leq\epsilon\frac{1-\exp(-T)}{1-\exp(-1)}
\end{multline*}
and
\begin{equation*}
    \left|f_d(x)\right|\leq\frac{\epsilon}{2\sqrt{\pi}}\int_{-\infty}^\infty \frac{1+\lambda^2}{1-\exp(-1)}\exp\left(-\frac{\lambda^2}{4}\right)d\lambda = \frac{3\epsilon}{1-\exp(-1)}.
\end{equation*}
After a simple calculation, we see that the integrals \eqref{int1} and \eqref{int2} satisfy the equation and the conditions (perturbation problem). Indeed, integrals \eqref{int1} and \eqref{int2} can be represented as follows
\begin{multline*}
    u_d(t,x)=\frac{\epsilon}{2\sqrt{\pi}}\int_{-\infty}^\infty \frac{1-\exp(-t)}{1-\exp(-1)}\exp\left(-\frac{\lambda^2}{4}\right)\exp(ix\lambda)d\lambda\\
    = \frac{\epsilon}{2\sqrt{\pi}}\frac{1-\exp(-t)}{1-\exp(-1)}\int_{-\infty}^\infty\exp\left(-\frac{\lambda^2}{4}\right)\exp(ix\lambda)d\lambda = \frac{1-\exp(-t)}{1-\exp(-1)}\epsilon\cdot\exp(-x^2)
\end{multline*}
and
\begin{multline*}
    f_d(x)=\frac{\epsilon}{2\sqrt{\pi}}\int_{-\infty}^\infty \frac{1+\lambda^2}{1-\exp(-1)}\exp\left(-\frac{\lambda^2}{4}\right)\exp(ix\lambda)d\lambda\\
    =\frac{\epsilon}{2\sqrt{\pi}}\frac{1}{1-\exp(-1)}\int_{-\infty}^\infty (1+\lambda^2)\exp\left(-\frac{\lambda^2}{4}\right)\exp(ix\lambda)d\lambda\\
    = \frac{1}{1-\exp(-1)}\epsilon\cdot\exp(-x^2)(3-4x^2).
\end{multline*}
In the following pictures you can find the graphics of the functions $f_d(x)= \frac{1}{1-\exp(-1)}\epsilon\cdot\exp(-x^2)(3-4x^2)$ and $u_d(x,y)=\frac{1-\exp(-x)}{1-\exp(-1)}\epsilon\cdot\exp(-y^2)$ for different epsilons ($\epsilon=1,0.5,0.1$). 

\begin{figure}[h]
    \centering
    \includegraphics[width=0.75\textwidth]{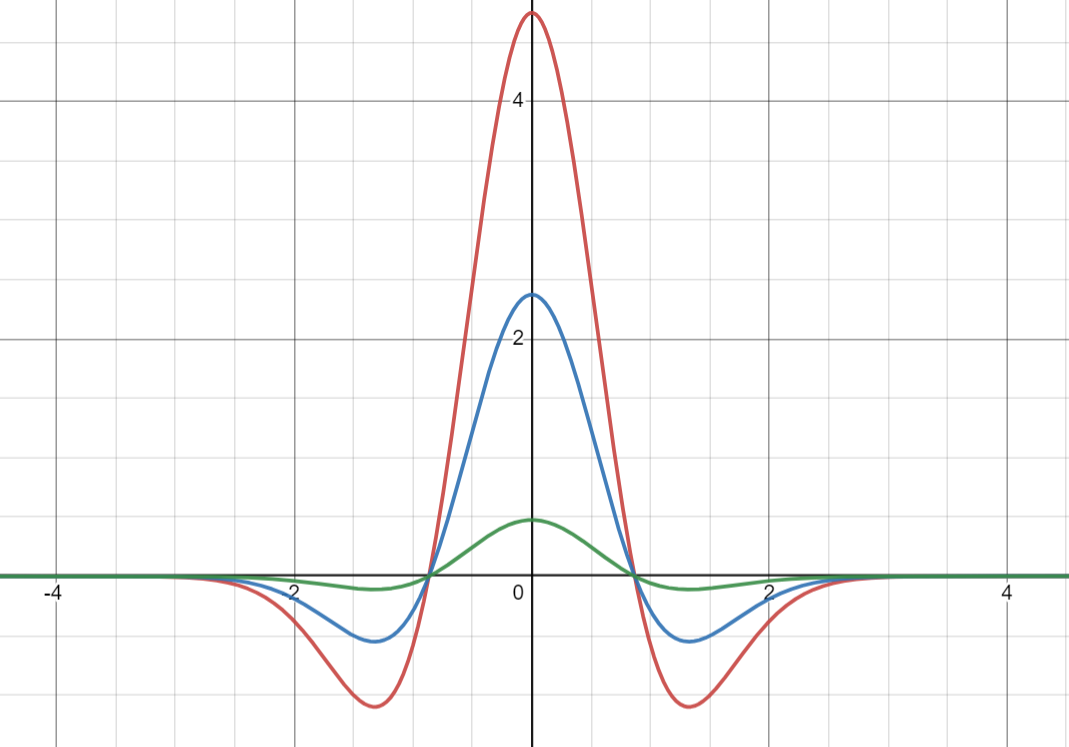}
    \caption{The graph of the function $f_d$, here we have used desmos.com. The red graph with $\epsilon=1$, the blue graph with $\epsilon=0.5$, and the green graph with $\epsilon=0.1$.}
    \label{fig:mesh1}
\end{figure}

\newpage

\begin{figure}[h]
    \centering
    \includegraphics[width=0.75\textwidth]{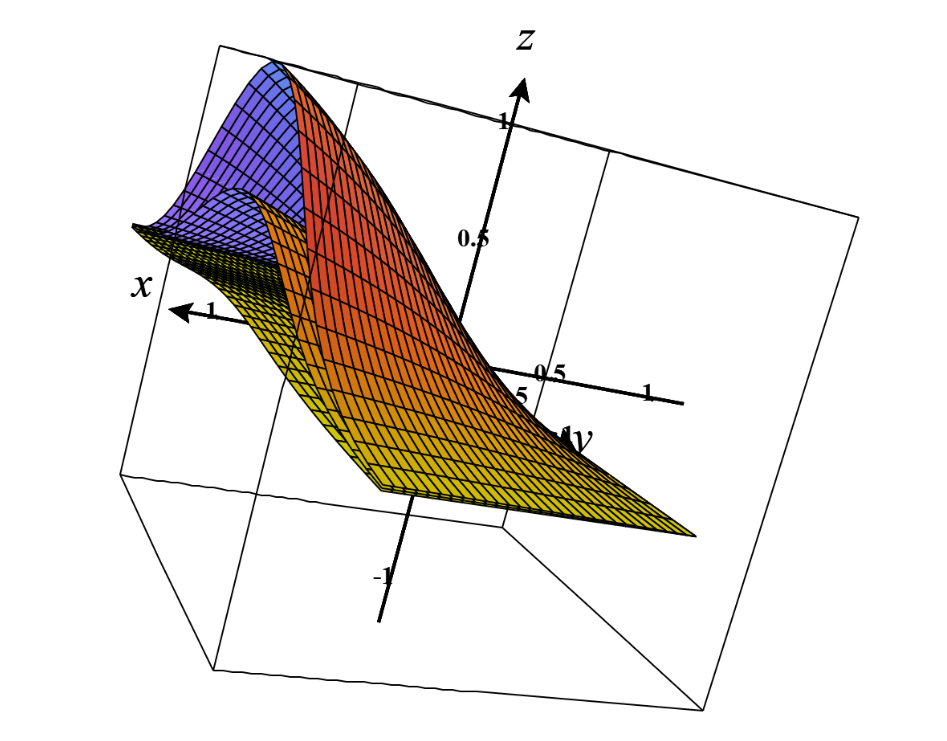}
    \caption{The graph of the function $u_d$, here we have used 3D Calc Plotter. The upper graph with $\epsilon=1$, the middle graph with $\epsilon=0.5$, and the lower graph with $\epsilon=0.1$.}
    \label{fig:mesh1}
\end{figure}

Now, let us calculate the following integrals:
\begin{multline*}
    \|\psi-\psi_d\|_{\mathcal{H}_\alpha^2} = \left(\frac{1}{\sqrt{2\pi}} \int_\mathbb{R} \left|(1+\lambda^2)\frac{\epsilon}{\sqrt{2}}\exp\left(-\frac{\lambda^2}{4}\right)\right|^2 d\lambda\right)^{\frac{1}{2}}\\
    = \epsilon\left(\frac{1}{2\sqrt{2\pi}} \int_\mathbb{R} (1+\lambda^2)^2\exp\left(-\frac{\lambda^2}{2}\right)d\lambda\right)^{\frac{1}{2}} = \epsilon\sqrt{3},
\end{multline*}
since $\mathcal{F}(\psi-\psi_d)(\lambda)=-\frac{\epsilon}{\sqrt{2}}\exp\left(-\frac{\lambda^2}{4}\right)$,
\begin{multline*}
    \|f-f_d\|_{2,\alpha} = \left(\frac{1}{\sqrt{2\pi}} \int_\mathbb{R}\left|f(x)-f_d(x)\right|^2 dx\right)^{\frac{1}{2}}\\
    = \epsilon\left(\frac{1}{\sqrt{2\pi}(1-\exp(-1))^2} \int_\mathbb{R}\left|\exp(-x^2)(3-4x^2)\right|^2dx\right)^{\frac{1}{2}} = \epsilon\frac{\sqrt{3}\exp(1)}{\exp(1)-1},
\end{multline*}
and 
\begin{equation*}
    \|u-u_d\|_{C([0,1],\mathcal{H}_\alpha^2(\mathbb{R},\mu_\alpha))} = \max_{0\leq t\leq 1} \|u(t,\cdot)-u_d(t,\cdot)\|_{\mathcal{H}_\alpha^2} = \epsilon\sqrt{3} \max_{0\leq t\leq 1}\frac{1-\exp(-t)}{1-\exp(-1)} = \epsilon\sqrt{3},
\end{equation*}
since 
\begin{multline*}
    \|u(t,\cdot)-u_d(t,\cdot)\|_{\mathcal{H}_\alpha^2} = \left(\frac{1}{\sqrt{2\pi}} \int_\mathbb{R}(1+\lambda^2)^2\left|\widehat{u}(t,\lambda)-\widehat{u}_d(t,\lambda)\right|^2 d\lambda\right)^{\frac{1}{2}}\\
    = \left(\frac{1}{\sqrt{2\pi}} \int_\mathbb{R}(1+\lambda^2)^2\left|\frac{\epsilon}{\sqrt{2}}\frac{1-\exp(-t)}{1-\exp(-1)}\exp\left(-\frac{\lambda^2}{4}\right)\right|^2 d\lambda\right)^{\frac{1}{2}}\\
    = \left(\frac{\epsilon^2}{2\sqrt{2\pi}}\left(\frac{1-\exp(-t)}{1-\exp(-1)}\right)^2\int_\mathbb{R}(1+\lambda^2)^2\exp\left(-\frac{\lambda^2}{2}\right) d\lambda\right)^{\frac{1}{2}}\\
    = \epsilon\frac{1-\exp(-t)}{1-\exp(-1)}\left(\frac{1}{2\sqrt{2\pi}}\int_\mathbb{R}(1+\lambda^2)^2\exp\left(-\frac{\lambda^2}{2}\right) d\lambda\right)^{\frac{1}{2}}.
\end{multline*}
According to the above computations we able to build a according table for different values of epsilon ($\epsilon=1,0.5,0.1$). 

\begin{table}[h]
    \centering
    \begin{tabular}{|p{4cm}|p{3cm}|p{3cm}|p{3cm}|}
      \hline
      $\epsilon$ & 1 & 0.5 & 0.1 \\
      \hline
      $\|\psi-\psi_d\|_{\mathcal{H}_\alpha^2}$ & 1.73205 & 0.86602 & 0.173205 \\ 
      \hline
      $\|f-f_d\|_{2,\alpha}$ & 2.74006 & 1.37003 & 0.274006 \\ 
      \hline
      $\|u-u_d\|_{C([0,1],\mathcal{H}_\alpha^2(\mathbb{R},\mu_\alpha))}$ & 1.73205 & 0.86602 & 0.173205 \\ 
      \hline
    \end{tabular}
    \caption{Stability test}
    \label{tab1}
\end{table}

\textbf{Conclusion}. In this subsection, we have considered one ISP and defined its solution using our calculus developed in this paper. Then to examine stability of its solution, we considered perturbation problem. The table \ref{tab1} shows that solution of the ISP is stable regarding to the small changes of the dates.

\end{document}